\documentclass[11pt]{article}
\usepackage{amssymb}
\usepackage{graphicx}
\title{Formulas and properties\footnote{The study was carried out within the framework of the state
contract of the Sobolev Institute of Mathematics (project No.
0314-2019-0002), and the Committee of Science in Education and the
Science Ministry of the Republic of Kazakhstan (Grant No.
AP08855544).}}
\author{S.V. Sudoplatov}
\date{}

\begin{document}

\maketitle
\begin{abstract}
We study links between first-order formulas and arbitrary
properties for families of theories, classes of structures and
their isomorphism types. Possibilities for ranks and degrees for
formulas and theories with respect to given properties are
described. Characteristics for generic sentences and generic
theories with respect to properties are described and
characterized.
\end{abstract}


\bigskip

First-order formulas are used to express semantic and syntactic
definable properties. Since in general there are more natural
properties than definable ones, these formulas can express them in
a partial way. In the present paper we study links between
formulas and arbitrary properties, consider characteristics
reflecting measures of their correspondence.

The paper is organized as follows. Preliminary notions, notations
and results are represented in Section 1. In Section 2 we consider
links between formulas and properties for semantic and syntactic
families. We study and characterize these links with respect to
constructions of formulas, set-theoretic operations and closures.
In Section 3 we study, characterize and describe rank values and
degree values for formulas with respect to given properties. In
Section 4, we describe spectra for cardinalities of definable
properties. In Section 5, generic sentences and theories with
respect to properties are introduced, and their links and ranks
and described.

Throughout we use the standard terminology in Mathematical Logic
\cite{ErPa}, notions, notations and results of \cite{RSrank,
MarkhSud, dst}.

\section{Preliminaries}

Let $\Sigma$ be a language. If $\Sigma$ is relational we denote by
$\mathcal{T}_\Sigma$ the family of all theories of the language
$\Sigma$. If $\Sigma$ contains functional symbols $f$ then
$\mathcal{T}_\Sigma$ is the family of all theories of the language
$\Sigma'$, which is obtained by replacements of all $n$-ary
symbols $f$ with $(n+1)$-ary predicate symbols $R_f$ interpreted
by $R_f=\{(\bar{a},b)\mid f(\bar{a})=b\}$.

Following \cite{RSrank} we define the {\em rank} ${\rm RS}(\cdot)$
for families $\mathcal{T}\subseteq\mathcal{T}_\Sigma$, similar to
Morley rank for a fixed theory, and a hierarchy with respect to
these ranks in the following way.

By $F(\Sigma)$ we denote the set of all formulas in the language
$\Sigma$ and by ${\rm Sent}(\Sigma)$ the set of all sentences in
$F(\Sigma)$.

For a sentence $\varphi\in {\rm Sent}(\Sigma)$ we denote by
$\mathcal{T}_\varphi$ the set of all theories $T\in\mathcal{T}$
with $\varphi\in T$.

Any set $\mathcal{T}_\varphi$ is called the {\em
$\varphi$-neighbourhood}, or simply a {\em neighbourhood}, for
$\mathcal{T}$, or the ($\varphi$-){\em definable} subset of
$\mathcal{T}$. The set $\mathcal{T}_\varphi$ is also called ({\em
formula-} or {\em sentence-}){\em definable} (by the sentence
$\varphi$) with respect to $\mathcal{T}$, or ({\em
sentence-})$\mathcal{T}$-{\em definable}, or simply {\em
$s$-definable}.

\medskip
{\bf Definition} \cite{RSrank}. For the empty family $\mathcal{T}$
we put the rank ${\rm RS}(\mathcal{T})=-1$, for finite nonempty
families $\mathcal{T}$ we put ${\rm RS}(\mathcal{T})=0$, and for
infinite families $\mathcal{T}$~--- ${\rm RS}(\mathcal{T})\geq 1$.

For a family $\mathcal{T}$ and an ordinal $\alpha=\beta+1$ we put
${\rm RS}(\mathcal{T})\geq\alpha$ if there are pairwise
inconsistent $\Sigma(\mathcal{T})$-sentences $\varphi_n$,
$n\in\omega$, such that ${\rm
RS}(\mathcal{T}_{\varphi_n})\geq\beta$, $n\in\omega$.

If $\alpha$ is a limit ordinal then ${\rm
RS}(\mathcal{T})\geq\alpha$ if ${\rm RS}(\mathcal{T})\geq\beta$
for any $\beta<\alpha$.

We set ${\rm RS}(\mathcal{T})=\alpha$ if ${\rm
RS}(\mathcal{T})\geq\alpha$ and ${\rm
RS}(\mathcal{T})\not\geq\alpha+1$.

If ${\rm RS}(\mathcal{T})\geq\alpha$ for any $\alpha$, we put
${\rm RS}(\mathcal{T})=\infty$.

A family $\mathcal{T}$ is called {\em $e$-totally transcendental},
or {\em totally transcendental}, if ${\rm RS}(\mathcal{T})$ is an
ordinal.

If $\mathcal{T}$ is $e$-totally transcendental, with ${\rm
RS}(\mathcal{T})=\alpha\geq 0$, we define the {\em degree} ${\rm
ds}(\mathcal{T})$ of $\mathcal{T}$ as the maximal number of
pairwise inconsistent sentences $\varphi_i$ such that ${\rm
RS}(\mathcal{T}_{\varphi_i})=\alpha$.

\medskip
{\bf Definition} \cite{at}. An infinite family $\mathcal{T}$ is
called {\em $e$-minimal} if for any sentence
$\varphi\in\Sigma(\mathcal{T})$, $\mathcal{T}_\varphi$ is finite
or $\mathcal{T}_{\neg\varphi}$ is finite.

\medskip
By the definition a family $\mathcal{T}$ is $e$-minimal iff ${\rm
RS}(\mathcal{T})=1$ and ${\rm ds}(\mathcal{T})=1$ \cite{RSrank},
and iff $\mathcal{T}$ has unique accumulation point \cite{at}.

\medskip
In the paper \cite{cl} the notion of $E$-closure was introduced
and characterized as follows:

\medskip
{\bf Proposition 1.1.} {\em If
$\mathcal{T}\subseteq\mathcal{T}_\Sigma$ is an infinite set and
$T\in\mathcal{T}_\Sigma\setminus\mathcal{T}$ then $T\in{\rm
Cl}_E(\mathcal{T})$ {\rm (}i.e., $T$ is an {\sl accumulation
point} for $\mathcal{T}$ with respect to $E$-closure ${\rm
Cl}_E${\rm )} if and only if for any sentence $\varphi\in T$ the
set $\mathcal{T}_\varphi$ is infinite.}

\medskip

The following theorem characterizes the property of $e$-total
transcendency for countable languages.

\medskip
{\bf Theorem 1.2} \cite{RSrank}. {\em For any family $\mathcal{T}$
with $|\Sigma(\mathcal{T})|\leq\omega$ the following conditions
are equivalent:

$(1)$ $|{\rm Cl}_E(\mathcal{T})|=2^\omega$;

$(2)$ $e$-${\rm Sp}(\mathcal{T})=2^\omega$;

$(3)$ ${\rm RS}(\mathcal{T})=\infty$.}

\medskip
{\bf Theorem 1.3} \cite{MarkhSud}. {\em For any language $\Sigma$
either ${\rm RS}(\mathcal{T}_\Sigma)$ is finite, if $\Sigma$
consists of finitely many $0$-ary and unary predicates, and
finitely many constant symbols, or ${\rm
RS}(\mathcal{T}_\Sigma)=\infty$, otherwise.}

\medskip
For a language $\Sigma$ we denote by $\mathcal{T}_{\Sigma,n}$ the
family of all theories in $\mathcal{T}_\Sigma$ having $n$-element
models, $n\in\omega$, as well as by $\mathcal{T}_{\Sigma,\infty}$
the family of all theories in $\mathcal{T}_\Sigma$ having infinite
models.

\medskip
{\bf Theorem 1.4} \cite{MarkhSud}. {\em For any language $\Sigma$
either ${\rm RS}(\mathcal{T}_{\Sigma,n})=0$, if $\Sigma$ is finite
or $n=1$ and $\Sigma$ has finitely many predicate symbols, or
${\rm RS}(\mathcal{T}_{\Sigma,n})=\infty$, otherwise.}

\medskip
{\bf Theorem 1.5} \cite{MarkhSud}. {\em For any language $\Sigma$
either ${\rm RS}(\mathcal{T}_{\Sigma,\infty})$ is finite, if
$\Sigma$ is finite and without predicate symbols of arities $m\geq
2$ as well as without functional symbols of arities $n\geq 1$, or
${\rm RS}(\mathcal{T}_{\Sigma,\infty})=\infty$, otherwise.}

\medskip
By the definition the families $\mathcal{T}_\Sigma$,
$\mathcal{T}_{\Sigma,n}$, $\mathcal{T}_{\Sigma,\infty}$ are
$E$-closed. Thus, combining Theorem 1.2 with Theorems 1.3--1.5 we
obtain the following possibilities of cardinalities for the
families $\mathcal{T}_\Sigma$, $\mathcal{T}_{\Sigma,n}$,
$\mathcal{T}_{\Sigma,\infty}$ depending on $\Sigma$ and
$n\in\omega$:

\medskip
{\bf Proposition 1.6.} {\em For any language $\Sigma$ either
either $\mathcal{T}_\Sigma$ is countable, if if $\Sigma$ consists
of finitely many $0$-ary and unary predicates, and finitely many
constant symbols, or $|\mathcal{T}_\Sigma|\geq 2^\omega$,
otherwise.}

\medskip
{\bf Proposition 1.7.} {\em For any language $\Sigma$ either
$\mathcal{T}_{\Sigma,n}$ is finite, if $\Sigma$ is finite or $n=1$
and $\Sigma$ has finitely many predicate symbols, or
$|\mathcal{T}_{\Sigma,n}|\geq 2^\omega$, otherwise.}

\medskip
{\bf Proposition 1.8.} {\em For any language $\Sigma$ either
$\mathcal{T}_{\Sigma,\infty}$ is at most countable, if $\Sigma$ is
finite and without predicate symbols of arities $m\geq 2$ as well
as without functional symbols of arities $n\geq 1$, or
$|\mathcal{T}_{\Sigma,\infty}|\geq 2^\omega$, otherwise.}

\medskip
{\bf Definition} \cite{dst}. If $\mathcal{T}$ is a family of
theories and $\Phi$ is a set of sentences, then we put
$\mathcal{T}_\Phi=\bigcap\limits_{\varphi\in\Phi}\mathcal{T}_\varphi$
and the set $\mathcal{T}_\Phi$ is called ({\em type-}) or ({\em
diagram-}){\em definable} (by the set $\Phi$) with respect to
$\mathcal{T}$, or ({\em diagram-})$\mathcal{T}$-{\em definable},
or simply {\em $d$-definable}.

Clearly, finite unions of $d$-definable sets are again
$d$-definable. Considering infinite unions $\mathcal{T}'$ of
$d$-definable sets $\mathcal{T}_{\Phi_i}$, $i\in I$, one can
represent them by sets of sentences with infinite disjunctions
$\bigvee\limits_{i\in I}\varphi_i$, $\varphi_i\in\Phi_i$. We call
these unions $\mathcal{T}'$ are called {\em $d_\infty$-definable}
sets.

\medskip
{\bf Definition} \cite{dst}. Let $\mathcal{T}$ be a family of
theories, $\Phi$ be a set of sentences, $\alpha$ be an ordinal
$\leq{\rm RS}(\mathcal{T})$ or $-1$. The set $\Phi$ is called {\em
$\alpha$-ranking} for $\mathcal{T}$ if ${\rm
RS}(\mathcal{T}_\Phi)=\alpha$. A sentence $\varphi$ is called {\em
$\alpha$-ranking} for $\mathcal{T}$ if $\{\varphi\}$ is {\em
$\alpha$-ranking} for $\mathcal{T}$.

The set $\Phi$ (the sentence $\varphi$) is called {\em ranking}
for $\mathcal{T}$ if it is $\alpha$-ranking for $\mathcal{T}$ with
some $\alpha$.

\medskip
{\bf Proposition 1.9} \cite{dst}. {\em For any ordinals
$\alpha\leq\beta$, if ${\rm RS}(\mathcal{T})=\beta$ then ${\rm
RS}(\mathcal{T}_\varphi)=\alpha$ for some {\rm
(}$\alpha$-ranking{\rm )} sentence $\varphi$. Moreover, there are
${\rm ds}(\mathcal{T})$ pairwise $\mathcal{T}$-inconsistent
$\beta$-ranking sentences for $\mathcal{T}$, and if $\alpha<\beta$
then there are infinitely many pairwise $\mathcal{T}$-inconsistent
$\alpha$-ranking sentences for $\mathcal{T}$.}

\medskip
{\bf Theorem 1.10} \cite{dst}. {\em Let $\mathcal{T}$ be a family
of a countable language $\Sigma$ and with ${\rm
RS}(\mathcal{T})=\infty$, $\alpha$ be a countable ordinal,
$n\in\omega\setminus\{0\}$. Then there is a $d_\infty$-definable
subfamily $\mathcal{T}^\ast\subset\mathcal{T}$ such that ${\rm
RS}(\mathcal{T}^\ast)=\alpha$ and ${\rm ds}(\mathcal{T}^\ast)=n$.}

\medskip
{\bf Theorem 1.11} \cite{hataf}. {\em For any two disjoint
subfamilies $\mathcal{T}_1$ and $\mathcal{T}_2$ of an $E$-closed
family $\mathcal{T}$ the following conditions are equivalent:

$(1)$ $\mathcal{T}_1$ and $\mathcal{T}_2$ are separated by some
sentence $\varphi$: $\mathcal{T}_1\subseteq\mathcal{T}_\varphi$
and $\mathcal{T}_2\subseteq\mathcal{T}_{\neg\varphi}$;

$(2)$ $E$-closures of $\mathcal{T}_1$ and $\mathcal{T}_2$  are
disjoint in $\mathcal{T}$: ${\rm Cl}_E(\mathcal{T}_1)\cap{\rm
Cl}_E(\mathcal{T}_2)\cap\mathcal{T}=\emptyset$;

$(3)$ $E$-closures of $\mathcal{T}_1$ and $\mathcal{T}_2$ are
disjoint: ${\rm Cl}_E(\mathcal{T}_1)\cap{\rm
Cl}_E(\mathcal{T}_2)=\emptyset$.}

\medskip
{\bf Definition} \cite{cl}. Let $\mathcal{T}_0$ be a family of
theories. A subset $\mathcal{T}'_0\subseteq\mathcal{T}_0$ is said
to be {\em generating}\index{Set!generating} if
$\mathcal{T}_0={\rm Cl}_E(\mathcal{T}'_0)$. The generating set
$\mathcal{T}'_0$ (for $\mathcal{T}_0$) is {\em
minimal}\index{Set!generating!minimal} if $\mathcal{T}'_0$ does
not contain proper generating subsets. A minimal generating set
$\mathcal{T}'_0$ is {\em least}\index{Set!generating!least} if
$\mathcal{T}'_0$ is contained in each generating set for
$\mathcal{T}_0$.

\medskip
{\bf Theorem 1.12} \cite{cl}. {\em If $\mathcal{T}'_0$ is a
generating set for a $E$-closed set $\mathcal{T}_0$ then the
following conditions are equivalent:

$(1)$ $\mathcal{T}'_0$ is the least generating set for
$\mathcal{T}_0$;

$(2)$ $\mathcal{T}'_0$ is a minimal generating set for
$\mathcal{T}_0$;

$(3)$ any theory in $\mathcal{T}'_0$ is isolated by some set
$(\mathcal{T}'_0)_\varphi$, i.e., for any $T\in\mathcal{T}'_0$
there is $\varphi\in T$ such that
$(\mathcal{T}'_0)_\varphi=\{T\}$;

$(4)$ any theory in $\mathcal{T}'_0$ is isolated by some set
$(\mathcal{T}_0)_\varphi$, i.e., for any $T\in\mathcal{T}'_0$
there is $\varphi\in T$ such that
$(\mathcal{T}_0)_\varphi=\{T\}$.}

\section{Relations between formulas and properties}

\medskip
{\bf Definition.} Let $\Sigma$ be a language,
$\varphi\rightleftharpoons\varphi(\overline{x})$ be a formula in
$F(\Sigma)$, $P_s$ be a subclass of the class $K(\Sigma)$ of all
structures $\mathcal{A}$ in the language $\Sigma$. We say that
$\varphi(\overline{x})$ {\em partially} (respectively, {\em
totally}) {\em satisfies} $P_s$, denoted by $\varphi\rhd_{\rm ps}
P_s$ or $\varphi\rhd^\exists_s P_s$ ($\varphi\rhd_{\rm ts} P_s$ or
$\varphi\rhd^\forall_s P_s$), if there are $\mathcal{A}\in P_s$
and $\overline{a}\in A$ (for any $\mathcal{A}\in P_s$ there is
$\overline{a}\in A$) such that
$\mathcal{A}\models\varphi(\overline{a})$.

If $P_{\rm is}$ is a subclass of the class ${\rm ITK}(\Sigma)$ of
isomorphism types for the class $K(\Sigma)$ then we say that
$\varphi(\overline{x})$ {\em partially} (respectively, {\em
totally}) {\em satisfies} $P_{\rm its}$, denoted by
$\varphi\rhd_{\rm pits} P_{\rm its}$ or $\varphi\rhd^\exists_{\rm
its} P_{\rm its}$ ($\varphi\rhd_{\rm tits} P_{\rm its}$ or
$\varphi\rhd^\forall_{\rm its} P_{\rm its}$) if $\varphi\rhd_{\rm
ps} P_s$ ($\varphi\rhd_{\rm ts} P_s$, where $P_s$ consists of all
structures whose isomorphism types belong to $P_{\rm its}$.

If $P_t$ is a subset of the set $\mathcal{T}_\Sigma$ of all
complete theories in the language $\Sigma$ then we say that
$\varphi(\overline{x})$ {\em partially} (respectively, {\em
totally}) {\em satisfies} $P_t$, denoted by $\varphi\rhd_{\rm pt}
P_t$ or $\varphi\rhd^\exists_t P_t$ ($\varphi\rhd_{\rm tt} P_t$ or
$\varphi\rhd^\forall_t P_t$), if there are $T\in P_t$,
$\mathcal{M}\models T$, and $\overline{a}\in M$ (for any $T\in
P_t$ there are $\mathcal{M}\models T$ and $\overline{a}\in M$)
such that $\mathcal{M}\models\varphi(\overline{a})$.

We write $\not\!\rhd_\xi$ if a $\rhd_\xi$-relation does not hold.

\medskip
{\bf Remark 2.1.} By the definition we have the following obvious
properties.

1. If $P_s\ne\emptyset$ and $\varphi\rhd_{\rm ts} P_s$ then
$\varphi\rhd_{\rm ps} P_s$. Similarly, if $\varphi\rhd_{\rm tits}
P_{\rm its}$ for nonempty $P_{\rm its}$ then $\varphi\rhd_{\rm
pits} P_{\rm its}$, and if $\varphi\rhd_{\rm tt} P_t$ for nonempty
$P$ then $\varphi\rhd_{\rm pt} P_t$.

\medskip
2. For any singleton $P_s$, $\varphi\rhd_{\rm ps} P_s$ implies
$\varphi\rhd_{\rm ts} P_s$. Similarly, $\varphi\rhd_{\rm pits}
P_{\rm its}$ implies $\varphi\rhd_{\rm tits} P_{\rm its}$ for any
singleton $P_{\rm its}$, and $\varphi\rhd_{\rm pt} P_t$ implies
$\varphi\rhd_{\rm tt} P$ for any singleton $P_t$.

\medskip
3. If $\varphi\rhd_{\rm ps} \{\mathcal{A}\}$ and
$\mathcal{A}\equiv\mathcal{B}$ then $\varphi\rhd_{\rm ps}
\{\mathcal{B}\}$. It implies that the relations $\varphi\rhd_{\rm
pt} P_t$ and $\varphi\rhd_{\rm tt} P_t$ do not depend on the
choice of models $\mathcal{M}\models T$ for $T\in P_t$.

\medskip
4. (Reflexivity) For any sentence $\varphi$ and a (nonempty)
family $\mathcal{T}_\varphi\subseteq\mathcal{T}_\Sigma$ we have
$\varphi\rhd_{\rm tt} \mathcal{T}_\varphi$ (and $\varphi\rhd_{\rm
pt} \mathcal{T}_\varphi$).

\medskip
5. (Monotony) If $\varphi\rhd_{\rm ps} P_s$, $\varphi\vdash\psi$
and $P_s\subseteq P'_s\subseteq K(\Sigma)$ then $\psi\rhd_{\rm ps}
P'_s$. If $\varphi\rhd_{\rm ts} P_s$, $\varphi\vdash\psi$ and
$P_s\supseteq P'_s$ then $\psi\rhd_{\rm ps} P'_s$. If
$\varphi\rhd_{\rm pits} P_{\rm its}$, $\varphi\vdash\psi$ and
$P_{\rm its}\subseteq P'_{\rm its}\subseteq {\rm ITK}(\Sigma)$
then $\psi\rhd_{\rm pits} P'_{\rm its}$. If $\varphi\rhd_{\rm
tits} P_{\rm its}$, $\varphi\vdash\psi$ and $P_{\rm its}\supseteq
P'_{\rm its}$ then $\psi\rhd_{\rm ps} P'_{\rm its}$. If
$\varphi\rhd_{\rm pt} P_t$, $\varphi\vdash\psi$ and $P_t\subseteq
P'_t\subseteq \mathcal{T}_\Sigma$ then $\psi\rhd_{\rm ps} P'_t$.
If $\varphi\rhd_{\rm tt} P_t$, $\varphi\vdash\psi$ and
$P_t\supseteq P'_t$ then $\psi\rhd_{\rm tt} P'_t$.

\medskip

For a property $P_s$ we denote by ${\rm ITK}(P_s)$ the class of
isomorphism types for structures in $P_s$, and by ${\rm Th}(P_s)$
the set $\{T\in \mathcal{T}_\Sigma\mid\mathcal{A}\models T$ for
some $\mathcal{A}\in P_s\}$.

For a property $P_{\rm its}$ we denote by $K(P_{\rm its})$ the
class of all structures whose isomorphism types are represented in
$P_{\rm its}$, and by ${\rm Th}(P_{\rm its})$ the set ${\rm
Th}(K(P_{\rm its}))$.

For a property $P_t$ we denote by $K(P_t)$ the class of all models
of theories in $P_t$, and by ${\rm ITK}(P_t)$ the class ${\rm
ITK}(K(P_t))$.

In terms of these notations by the definition we have the
following natural links between semantic properties $P_s$ and
$P_{\rm its}$ and syntactic properties $P_t$:

\medskip
{\bf Proposition 2.2.} {\em For any formula $\varphi\in F(\Sigma)$
and properties $P_s$, $P_{\rm its}$, $P_t$ the following
conditions hold:

$(1)$ $\varphi\rhd_{\rm ps} P_s$ iff $\varphi\rhd_{\rm pits} {\rm
ITK}(P_s)$, and iff $\varphi\rhd_{\rm pt}{\rm Th}(P_s)$;

$(2)$ $\varphi\rhd_{\rm ts} P_s$ iff $\varphi\rhd_{\rm tits} {\rm
ITK}(P_s)$, and iff $\varphi\rhd_{\rm tt}{\rm Th}(P_s)$;

$(3)$ $\varphi\rhd_{\rm pits} P_{\rm its}$ iff $\varphi\rhd_{\rm
ps} K(P_{\rm its})$, and iff $\varphi\rhd_{\rm pt}{\rm Th}(P_{\rm
its})$;

$(4)$ $\varphi\rhd_{\rm tits} P_{\rm its}$ iff $\varphi\rhd_{\rm
ts} K(P_{\rm its})$, and iff $\varphi\rhd_{\rm tt}{\rm Th}(P_{\rm
its})$;

$(5)$ $\varphi\rhd_{\rm pt} P_t$ iff $\varphi\rhd_{\rm ps}K(P_t)$,
and iff $\varphi\rhd_{\rm pits} {\rm ITK}(P_t)$;

$(6)$ $\varphi\rhd_{\rm tt} P_t$ iff $\varphi\rhd_{\rm ts}K(P_t)$,
and iff $\varphi\rhd_{\rm tits} {\rm ITK}(P_t)$.

In the items $(3)$ and $(4)$ the class $K(P_{\rm its})$ can be
replaced by a subclass $K'$ such that ${\rm ITK}(K')=P_{\rm its}$.
Similarly, in the items $(5)$ and $(6)$ the class $K(P_t)$ can be
replaced by a subclass $K'$ such that ${\rm Th}(K')=P_t$, and
independently ${\rm ITK}(P_t)$ can be replaced by a subclass $K''$
such that ${\rm Th}(K'')=P_t$.}

\medskip
By Proposition 2.2 semantic properties $P_s$ and $P_{\rm its}$ can
be naturally transformed into syntactic ones $P_t$, and vice
versa. It means that natural model-theoretic properties such as
$\omega$-categoricity, stability, simplicity etc. can be
formulated both for theories, for structures and for their
isomorphism types.

The links between $\rhd$-relations which pointed out in
Proposition 2.2 allow to reduce our consideration to the relations
$\rhd_{\rm pt}$ and $\rhd_{\rm tt}$. Besides, for the simplicity
we will principally consider sentences $\varphi$ instead of
formulas in general. Reductions of formulas $\psi(\overline{x})$
to sentences use the operators
$\psi(\overline{x})\mapsto\forall\overline{x}\,\psi(\overline{x})$
and
$\psi(\overline{x})\mapsto\exists\overline{x}\,\psi(\overline{x})$.

\medskip
{\bf Proposition 2.3.} {\em For any sentences $\varphi,\psi\in
{\rm Sent}(\Sigma)$ and properties $P_t,P'_t\subseteq
\mathcal{T}_\Sigma$ the following conditions hold:

$(1)$ if $(\varphi\wedge\psi)\rhd_{\rm pt} (P_t\cap P'_t)$ then
$\varphi\rhd_{\rm pt} P_t$ and $\psi\rhd_{\rm pt} P'_t$; the
converse implication does not hold: there are $\varphi',\psi'\in
{\rm Sent}(\Sigma)$ and $P''_t\in \mathcal{T}_\Sigma$ such that
$\varphi'\rhd_{\rm pt} P''_t$, $\psi'\rhd_{\rm pt} P''_t$, and
$(\varphi'\wedge\psi')\not\!\rhd_{\rm pt} P''_t$;

$(2)$ $(\varphi\wedge\psi)\rhd_{\rm tt} P_t$ iff $\varphi\rhd_{\rm
tt} P_t$ and $\psi\rhd_{\rm tt} P_t$;

$(3)$ if $\varphi\rhd_{\rm tt} P_t$ and $\psi\rhd_{\rm tt} P'_t$
then $(\varphi\wedge\psi)\rhd_{\rm tt} (P_t\cap P'_t)$; the
converse implication does not hold: there are $\varphi',\psi'\in
{\rm Sent}(\Sigma)$ and $P''_t,P'''_t\in \mathcal{T}_\Sigma$ such
that $(\varphi\wedge\psi)\rhd_{\rm tt} (P''_t\cap P'''_t)$ whereas
$\varphi'\not\!\rhd_{\rm tt} P''_t$ and $\psi'\not\!\rhd_{\rm tt}
P'''_t$.}

\medskip
Proof. $(1)$ If $(\varphi\wedge\psi)\rhd_{\rm pt} (P_t\cap P'_t)$
then there is $T\in P_t\cap P'_t$ with $(\varphi\wedge\psi)\in T$.
Since $\varphi,\psi\in T$, $T\in P_t$ and $T\in P'_t$, we obtain
$\varphi\rhd_{\rm pt} P_t$ and $\psi\rhd_{\rm pt} P'_t$. Therefore
it suffices to notice for $(1)$ that sentences $\varphi',\psi'$
asserting distinct finite cardinalities $m$ and $n$ for universes
partially satisfy a property $P''_t$ containing a theory $T_1$
with an $m$-element model and a theory $T_2$ with an $n$-element
model. At the same time $(\varphi'\wedge\psi')\not\!\rhd_{\rm pt}
P''_t$ since $(\varphi'\wedge\psi')$ is inconsistent.

$(2)$ If $(\varphi\wedge\psi)\rhd_{\rm tt} P_t$ then
$(\varphi\wedge\psi)$, and so, $\varphi$ and $\psi$ belong to all
theories in $P_t$ implying $\varphi\rhd_{\rm tt} P_t$ and
$\psi\rhd_{\rm tt} P_t$. Conversely, if $\varphi\rhd_{\rm tt} P_t$
and $\psi\rhd_{\rm tt} P_t$ then $\varphi$, $\psi$, and so
$(\varphi\wedge\psi)$ belong to all theories in $P_t$ implying
$(\varphi\wedge\psi)\rhd_{\rm tt} P_t$.

$(3)$ If $\varphi\rhd_{\rm tt} P_t$ and $\psi\rhd_{\rm tt} P'_t$
then $\varphi\in\bigcap P_t$ and $\psi\in\bigcap P'_t$, implying
$(\varphi\wedge\psi)\in \bigcap P_t\cap\bigcap P'_t$, i.e.,
$(\varphi\wedge\psi)\rhd_{\rm tt} (P_t\cap P'_t)$. Finally, if
$P''_t$ and $P'''_t$ are nonempty with $P''_t\cap
P'''_t=\emptyset$, and $\varphi'$, $\psi'$ are inconsistent
sentences then $(\varphi\wedge\psi)\rhd_{\rm tt} (P''_t\cap
P'''_t)$ and $\varphi'\not\!\rhd_{\rm pt} P''_t$ and
$\psi'\not\!\rhd_{\rm pt} P'''_t$ implying
$\varphi'\not\!\rhd_{\rm tt} P''_t$ and $\psi'\not\!\rhd_{\rm tt}
P'''_t$.~$\Box$

\medskip
{\bf Proposition 2.4.} {\em For any sentences $\varphi,\psi\in
{\rm Sent}(\Sigma)$ and properties $P_t,P'_t\subseteq
\mathcal{T}_\Sigma$ the following conditions hold:

$(1)$ if $\varphi\rhd_{\rm pt} P_t$ or $\psi\rhd_{\rm pt} P'_t$
then $(\varphi\vee\psi)\rhd_{\rm pt} (P_t\cup P'_t)$; the converse
implication does not hold: there are $\varphi',\psi'\in {\rm
Sent}(\Sigma)$ and $P''_t\in \mathcal{T}_\Sigma$ such that
$(\varphi'\vee\psi')\rhd_{\rm pt} P''_t$, and
$\varphi'\not\!\rhd_{\rm pt} P''_t$ or $\psi'\not\!\rhd_{\rm pt}
P''_t$;

$(2)$ $(\varphi\vee\psi)\rhd_{\rm pt} P_t$ iff $\varphi\rhd_{\rm
pt} P_t$ or $\psi\rhd_{\rm pt} P_t$;

$(3)$ if $\varphi\rhd_{\rm tt} P_t$ and $\psi\rhd_{\rm tt} P'_t$
then $(\varphi\vee\psi)\rhd_{\rm tt} (P_t\cup P'_t)$; the converse
implication does not hold: there are $\varphi',\psi'\in {\rm
Sent}(\Sigma)$ and $P''_t\in \mathcal{T}_\Sigma$ such that
$(\varphi'\vee\psi')\rhd_{\rm tt} (P''_t)$ whereas
$\varphi'\not\!\rhd_{\rm tt} P''_t$ and $\psi'\not\!\rhd_{\rm tt}
P''_t$.}

\medskip
Proof. $(1)$ If $\varphi\rhd_{\rm pt} P_t$ or $\psi\rhd_{\rm pt}
P'_t$ then $\varphi\in T$ for some $T\in P_t$ or $\psi\in T'$ for
some $T'\in P'_t$. Therefore $T$ or $T'$ witness that
$(\varphi\vee\psi)\rhd_{\rm pt} (P_t\cup P'_t)$. If $\varphi'$ is
a tautology and $\psi'$ is an inconsistent sentence, then for any
nonempty $P''_t\subseteq\mathcal{T}_\Sigma$ we have
$(\varphi'\vee\psi')\rhd_{\rm pt} P''_t$, $\varphi'\rhd_{\rm pt}
P''_t$, and $\psi'\not\!\rhd_{\rm pt} P''_t$.

$(2)$ It holds since a sentence $(\varphi\vee\psi)$ belongs to a
complete theory $T$ if and only if $\varphi\in T$ or $\psi\in T$.

$(3)$ If $\varphi\rhd_{\rm tt} P_t$ and $\psi\rhd_{\rm tt} P'_t$
then $\varphi\in\bigcap P_t$ and $\psi\in\bigcap P'_t$. Therefore
$(\varphi\vee\psi)\in\bigcap P_t$ and $(\varphi\vee\psi)\in\bigcap
P'_t$ implying $(\varphi\vee\psi)\in\bigcap P_t\cap\bigcap P'_t$,
hence $(\varphi\vee\psi)\rhd_{\rm tt} (P_t\cup P'_t)$.

Now let $P''_t=\mathcal{T}_\Sigma$, $\varphi'$ be a sentence
belonging to some but not all theories in $P''_t$. For the
sentence $\psi'=\neg\varphi'$ we have
$(\varphi'\vee\psi')\rhd_{\rm tt} (P''_t)$ since
$(\varphi\vee\psi)$ is a tautology, $\varphi'\not\!\rhd_{\rm tt}
P''_t$ and $\psi'\not\!\rhd_{\rm tt} P''_t$ by the choice of
$\varphi'$.~$\Box$

\medskip
{\bf Proposition 2.5.} {\em For any sentence $\varphi\in {\rm
Sent}(\Sigma)$ and a property $P_t\subseteq \mathcal{T}_\Sigma$
the following conditions hold:

$(1)$ $\varphi\rhd_{\rm pt} P_t$ iff $\neg\varphi\not\!\rhd_{\rm
tt} P_t$;

$(2)$ $\varphi\rhd_{\rm tt} P_t$ iff $\neg\varphi\not\!\rhd_{\rm
pt} P_t$.}

\medskip
Proof. $(1)$ If $\varphi\rhd_{\rm pt} P_t$ then there is $T\in
P_t$ such that $\varphi\in T$. Since $T$ is complete then
$\neg\varphi\notin T$ implying $\neg\varphi\not\!\rhd_{\rm tt}
P_t$. Conversely if $\neg\varphi\not\!\rhd_{\rm tt} P_t$ then
$\neg\varphi$ does not belong to some theory $T\in P_t$. Since $T$
is complete then $\varphi\in T$ implying $\varphi\rhd_{\rm pt}
P_t$.

$(2)$ is immediately follows from $(1)$.~$\Box$

\medskip
{\bf Proposition 2.6.} {\em For any formula $\varphi\in F(\Sigma)$
and a property $P_t\subseteq \mathcal{T}_\Sigma$ the following
conditions hold:

$(1)$ if $\varphi=\forall x\psi$ and $\varphi\rhd_{\rm pt} P_t$
then $\psi\rhd_{\rm pt} P_t$;

$(2)$ if $\varphi=\forall x\psi$ and $\varphi\rhd_{\rm tt} P_t$
then $\psi\rhd_{\rm tt} P_t$;

$(3)$ if $\varphi=\exists x\psi$ and $\psi\rhd_{\rm pt} P_t$ then
$\psi\rhd_{\rm pt} P_t$;

$(3)$ if $\varphi=\exists x\psi$ and $\psi\rhd_{\rm tt} P_t$ then
$\psi\rhd_{\rm tt} P_t$. }

\medskip
Proof. $(1)$ Let $\varphi=\forall x\psi$ and $\varphi\rhd_{\rm pt}
P_t$. Then there are $T\in P_t$, $\mathcal{M}\models T$,
$\overline{a}\in M$ such that
$\mathcal{M}\models\varphi(\overline{a})$. It implies
$\mathcal{M}\models\forall x\psi(x,\overline{a})$, therefore there
is $b\in M$ with $\mathcal{M}\models x\psi(b,\overline{a})$, i.e.,
$\psi\rhd_{\rm pt}\{T\}$, hence $\psi\rhd_{\rm pt} P_t$.

$(2)$ We repeat arguments for $(1)$ replacing some $T\in P_t$ by
an arbitrary one.

$(3)$ Let $\varphi=\exists x\psi$ and $\psi\rhd_{\rm pt} P_t$.
Then there are $T\in P_t$, $\mathcal{M}\models T$,
$\overline{a},b\in M$ such that
$\mathcal{M}\models\psi(b,\overline{a})$. It implies
$\mathcal{M}\models\exists x\psi(x,\overline{a})$, i.e.,
$\varphi\rhd_{\rm pt}\{T\}$, hence $\psi\rhd_{\rm pt} P_t$.

$(4)$ As for $(2)$ we repeat arguments for $(3)$ replacing some
$T\in P_t$ by an arbitrary one.~$\Box$

\medskip
The following two theorems assert that the relations $\rhd_{\rm
pt}$ and $\rhd_{\rm tt}$ are preserved under $E$-closures.

\medskip
{\bf Theorem 2.7.} {\em For any sentence $\varphi\in {\rm
Sent}(\Sigma)$ and a property $P_t\subseteq \mathcal{T}_\Sigma$
the following conditions are equivalent:

$(1)$ $\varphi\rhd_{\rm pt} P_t$,

$(2)$ $\varphi\rhd_{\rm pt}{\rm Cl}_E(P_t)$,

$(3)$ $\varphi\rhd_{\rm pt}P'_t$ for any/some $P'_t$ with ${\rm
Cl}_E(P'_t)={\rm Cl}_E(P_t)$. }

\medskip
Proof. $(1)\Rightarrow(2)$ holds in view of $P_t\subseteq{\rm
Cl}_E(P_t)$ and the monotony of the relation $\rhd_{\rm pt}$.

$(2)\Rightarrow(3)$. It suffices to show that $\varphi\rhd_{\rm
pt}P'_t$ for any $P'_t$ with ${\rm Cl}_E(P'_t)={\rm Cl}_E(P_t)$.
Since $\varphi\rhd_{\rm pt}{\rm Cl}_E(P_t)$ there is a theory
$T\in {\rm Cl}_E(P_t)$ with $\varphi\in T$. If $T\in P't$ we have
$\varphi\rhd_{\rm pt}P'_t$. Otherwise $T$ is an accumulation point
of $P'_t$ implying, in view of Proposition 1.1, that $varphi$
belongs to infinitely many theories in $P'_t$. Therefore
$\varphi\rhd_{\rm pt}P'_t$.

$(3)\Rightarrow(1)$ is obvious.~$\Box$

\medskip
{\bf Theorem 2.8.} {\em For any sentence $\varphi\in {\rm
Sent}(\Sigma)$ and a property $P_t\subseteq \mathcal{T}_\Sigma$
the following conditions are equivalent:

$(1)$ $\varphi\rhd_{\rm tt} P_t$,

$(2)$ $\varphi\rhd_{\rm tt}{\rm Cl}_E(P_t)$,

$(3)$ $\varphi\rhd_{\rm tt}P'_t$ for any/some $P'_t$ with ${\rm
Cl}_E(P'_t)={\rm Cl}_E(P_t)$. }

\medskip
Proof. $(1)\Rightarrow(2)$. Let $\varphi\rhd_{\rm tt} P_t$. If
$P_t$ is finite then ${\rm Cl}_E(P_t)=P_t$ and we have
$\varphi\rhd_{\rm tt}{\rm Cl}_E(P_t)$. If $P_t$ is finite then by
Proposition 1.1, ${\rm Cl}_E(P_t)$ consists of theories in $P_t$
and of theories $T\in\mathcal{T}_\Sigma\setminus P_t$ such that
for any sentence $\psi\in T$ the set $(P_t)_\psi$ is infinite.
Since $\varphi\rhd_{\rm tt} P_t$, $\varphi$ belongs to each such
theory $T$. Thus, $\varphi\rhd_{\rm tt}{\rm Cl}_E(P_t)$.

$(2)\Rightarrow(1)$ and $(2)\Rightarrow(3)$ are obvious.

$(3)\Rightarrow(2)$ follows assuming $\varphi\rhd_{\rm tt}P'_t$
for any/some $P'_t$ with ${\rm Cl}_E(P'_t)={\rm Cl}_E(P_t)$
repeating the arguments for $(1)\Rightarrow(2)$.~$\Box$

\medskip
For a property $P_t\subseteq \mathcal{T}_\Sigma$ we denote by
$\nabla(P_t)$ the set of all sentences $\varphi\in {\rm
Sent}(\Sigma)$ with $\varphi\rhd_{\rm pt} P_t$, and by
$\triangle(P_t)$ the set of all sentences $\varphi\in {\rm
Sent}(\Sigma)$ with $\varphi\rhd_{\rm tt} P_t$.

By the definition $\nabla(\emptyset)=\emptyset$,
$\triangle(\emptyset)={\rm Sent}(\Sigma)$,
$\nabla(\mathcal{T}_\Sigma)$ consists of all consistent sentences
$\varphi\in {\rm Sent}(\Sigma)$, and
$\triangle(\mathcal{T}_\Sigma)$ consists of all tautologies
$\varphi\in {\rm Sent}(\Sigma)$.

\medskip
{\bf Proposition 2.9.} {\em For any property $P_t\subseteq
\mathcal{T}_\Sigma$ the following conditions hold:

$(1)$ $\nabla(P_t)=\bigcup P_t$;

$(2)$ $\nabla(P_t)$ is consistent iff $|P_t|\leq 1$, and
$\nabla(P_t)$ is a complete theory iff $P_t$ is a singleton;

$(3)$ $\triangle(P_t)=\bigcap P_t$;

$(4)$ $\triangle(P_t)$ is a consistent theory iff
$P_t\ne\emptyset$, and $\triangle(P_t)$ is complete iff $P_t$ is a
singleton;

$(5)$ for any $P_t\ne\emptyset$,
$\nabla(P_t)\supseteq\triangle(P_t)$, and
$\nabla(P_t)=\triangle(P_t)$ iff $P_t$ is a singleton.}

\medskip
Proof. $(1)$ If $\varphi\in\nabla(P_t)$ then $\varphi\rhd_{\rm pt}
P_t$ and $\varphi\in T$ for some $T\in P_t$ implying
$\nabla(P_t)\subseteq\bigcup P_t$. Conversely, if
$\varphi\in\bigcup P_t$ then $\varphi\in T$ for some $T\in P_t$
implying $\varphi\rhd_{\rm pt} P_t$ and therefore
$\varphi\in\nabla(P_t)$.

$(2)$ Since $\nabla(\emptyset)=\emptyset$ it is consistent. If
$P_t=\{T\}$ then $\nabla(P_t)=T$, i.e., $\nabla(P_t)$ is
consistent and complete. If $P_t$ contains two distinct theories
$T_1$ and $T_2$ then $T_1\cup T_2\subseteq \nabla(P_t)$ implying
that $\nabla(P_t)$ is inconsistent since there are sentences
$\psi$ such that $\psi\in T_1$ and $\neg\psi\in T_2$.

$(3)$ If $\varphi\in\triangle(P_t)$ then $\varphi\rhd_{\rm tt}
P_t$ and $\varphi\in T$ for any $T\in P_t$ implying
$\triangle(P_t)\subseteq\bigcap P_t$. Conversely, if
$\varphi\in\bigcap P_t$ then $\varphi\in T$ for any $T\in P_t$
implying $\varphi\rhd_{\rm tt} P_t$ and therefore
$\varphi\in\triangle(P_t)$.

$(4)$ Since $\triangle(\emptyset)={\rm Sent}(\Sigma)$ it is
inconsistent. If $P_t\ne\emptyset$ then $\triangle(P_t)=\bigcap
P_t$ by $(3)$ implying that $\triangle(P_t)$ is a consistent
theory as an intersection of complete theories. If
$\triangle(P_t)$ is complete then both $P_t$ is nonempty and does
not contain two distinct theories, i.e., $P_t$ is a singleton.
Conversely, if $P_t=\{T\}$ then $\triangle(P_t)=T$ which is a
complete theory.

$(5)$ If $P_t\ne\emptyset$ then by $(1)$ and $(3)$ we have
$\nabla(P_t)=\bigcup P_t\supseteq\bigcap P_t=\triangle(P_t)$. If
$P_t=\{T\}$ then $\triangle(P_t)=T=\nabla(P_t)$. If $T_1,T_2\in
P_t$ for some $T_1\ne T_2$ then $\nabla(P_t)=\bigcup P_t\supseteq
T_1\cup T_2\supsetneq T_1\cap T_2\supseteq\bigcap
P_t=\triangle(P_t)$.~$\Box$

\medskip
Theorems 1.11, 2.7, 2.8 and Proposition 2.9 immediately imply the
following corollary on the separability of properties with respect
to the relations $\rhd_{\rm pt}$ and $\rhd_{\rm tt}$.

\medskip
{\bf Corollary 2.10.} {\em $(1)$ For any properties
$P_1,P_2\subseteq\mathcal{T}_\Sigma$ the following conditions
hold:

$(1)$ there exists $\varphi\in{\rm Sent}(\Sigma)$ such that
$\varphi\rhd_{\rm pt} P_1$ and $\neg\varphi\rhd_{\rm pt} P_2$ iff
$P_1$ and $P_2$ are nonempty and $|P_1\cup P_2|\geq 2$; in
particular, there exists $\varphi\in{\rm Sent}(\Sigma)$ such that
$\varphi\rhd_{\rm pt} P_1$ and $\neg\varphi\rhd_{\rm pt} P_1$ iff
$|P_1|\geq 2$;

$(2)$ there exists $\varphi\in{\rm Sent}(\Sigma)$ such that
$\varphi\rhd_{\rm tt} P_1$ and $\neg\varphi\rhd_{\rm tt} P_2$ iff
${\rm Cl}_E(P_1)\cap{\rm Cl}_E(P_2)=\emptyset$.}

\medskip
{\bf Corollary 2.11.} {\em For any nonempty property $P_t\subseteq
\mathcal{T}_\Sigma$ the following conditions hold:

$(1)$ the set $\triangle(P_t)$ forms a filter
$\triangle(P_t)/\!\!\equiv$ on
$\{\equiv\!\!(\varphi)\mid\varphi\in{\rm Sent}(\Sigma)\}$ with
respect to $\vdash$;

$(2)$ the filter $\triangle(P_t)/\!\!\equiv$ is principal iff
$\bigcap P_t$ is forced by some its sentence, i.e., $\bigcap P_t$
is a finitely axiomatizable theory, which is incomplete for
$|P_t|\geq 2$;

$(3)$ the filter $\triangle(P_t)/\!\!\equiv$ is an ultrafilter iff
$P_t$ is a singleton.}

\medskip
Proof. $(1)$ holds by Monotony and Proposition 2.3, (2).

$(2)$ immediately follows by Proposition 2.9, (3).

$(3)$ is satisfied in view of Proposition 2.9, (4).~$\Box$

\section{Ranks of sentences and spectra with respect to properties}

In this section we introduce a measure of complexity for sentences
satisfying a property using the ${\rm RS}$-rank for families of
theories \cite{RSrank, MarkhSud, dst}.

\medskip
{\bf Definition.} For a sentence $\varphi\in{\rm Sent}(\Sigma)$
and a property $P=P_t\subseteq\mathcal{T}_\Sigma$ we put ${\rm
RS}_P(\varphi)={\rm RS}(P_\varphi)$, and ${\rm ds}_P(\varphi)={\rm
ds}(P_\varphi)$ if ${\rm ds}(P_\varphi)$ is defined.

If $P=\mathcal{T}_\Sigma$ then we omit $P$ and write ${\rm
RS}(\varphi)$, ${\rm ds}(\varphi)$ instead of ${\rm
RS}_P(\varphi)$ and ${\rm ds}_P(\varphi)$, respectively.

\medskip
Clearly, if $P\subseteq P'\subseteq\mathcal{T}_\Sigma$ and
$\varphi\in{\rm Sent}(\Sigma)$ then ${\rm RS}_P(\varphi)\leq {\rm
RS}_{P'}(\varphi)$, and if ${\rm RS}_P(\varphi)={\rm
RS}_{P'}(\varphi)\in{\rm Ord}$ then ${\rm ds}_P(\varphi)\leq {\rm
ds}_{P'}(\varphi)$. Besides, we have:

\medskip
{\bf Proposition 3.1.} {\em $(1)$ $\varphi\rhd_{\rm tt} P$ iff
${\rm RS}_P(\neg\varphi)=-1$.

$(2)$ $\varphi\rhd_{\rm pt} P$ iff ${\rm RS}_P(\varphi)\geq 0$.}

\medskip
{\bf Definition.} For a sentence $\varphi\in{\rm Sent}(\Sigma)$
and a property $P\subseteq\mathcal{T}_\Sigma$ we say that
$\varphi$ is {\em $P$-totally transcendental} if ${\rm
RS}_P(\varphi)$ is an ordinal. If $P=\mathcal{T}_\Sigma$ then a
$P$-totally transcendental sentence is said to be {\em totally
transcendental}.

A sentence $\varphi$ is {\em co-{\rm (}$P${\rm )}-totally
transcendental} if $\neg\varphi$ is $P$-totally transcendental.

We omit $P$ and say about totally transcendental and co-totally
transcendental sentences if $P=\mathcal{T}_\Sigma$.

\medskip
By the definition each sentence $\varphi\in{\rm Sent}(\Sigma)$
obtains the characteristics ${\rm RS}_P(\varphi)$ and ${\rm
RS}_P(\neg\varphi)$ saying that $\varphi$ is (co)-rich enough with
respect to the property $P$, The characteristics ${\rm
ds}_P(\varphi)$ and ${\rm ds}_P(\neg\varphi)$, if they are
defined, give an additional information on ``$P$-richness'' of
$\varphi$.

For instance, if ${\rm RS}_P(\varphi)=0$ and ${\rm
ds}_P(\varphi)=n$ then $\varphi$ is {\em $P$-finite} satisfying
exactly $n$ theories in $P$. Respectively, if ${\rm
RS}_P(\neg\varphi)=0$ and ${\rm ds}_P(\neg\varphi)=n$ then
$\varphi$ is {\em $P$-cofinite}, i.e., it does not satisfy exactly
$n$ theories in $P$.

Clearly, $\varphi$ is both $P$-finite and $P$-cofinite iff $P$ is
nonempty and is finite.

\medskip
{\bf Theorem 3.2.} {\em For a language $\Sigma$ there is a totally
transcendental sentence $\varphi\in{\rm Sent}(\Sigma)$ iff
$\Sigma$ has finitely many predicate symbols.}

\medskip
Proof. If $\Sigma$ has finitely many predicate symbols we choose a
sentence $\varphi$ assigning that a universe is a singleton. Since
functional symbols have unique interpretations and there are
finitely many possibilities for (non)empty language predicate we
obtain ${\rm RS}(\varphi)=0$ that is $\varphi$ is totally
transcendental.

Conversely, if $\Sigma$ has infinitely many predicate symbols then
each consistent sentence $\varphi$ obtains a $2$-tree with respect
to (non)empty predicates $Q\in\Sigma\setminus\Sigma(\varphi)$.
This $2$-tree witnesses that ${\rm RS}(\varphi)=\infty$, i.e.,
$\varphi$ is not totally transcendental.~$\Box$.

\medskip
{\bf Remark 3.3.} If $\Sigma$ is finite then for the proof of
Theorem 3.2 it suffices to choose a sentence $\varphi$ assigning
that a universe is finite, since there are finitely many many
possibilities, up to isomorphism, for interpretations of language
symbols implying ${\rm RS}(\varphi)=0$.~$\Box$

\medskip
The following definition introduces values for richness of a
sentence with respect to a property.

\medskip
{\bf Definition.} For a language $\Sigma$, a property
$P\subseteq\mathcal{T}_\Sigma$, an ordinal $\alpha$ and a natural
number $n\geq 1$, a sentence $\varphi\in{\rm Sent}(\Sigma)$ is
called {\em $(P,\alpha,n)$-{\rm (}co-{\rm )}rich} if ${\rm
RS}_P(\varphi)=\alpha$ and ${\rm ds}_P(\varphi)=n$ (respectively,
${\rm RS}_P(\neg\varphi)=\alpha$ and ${\rm ds}_P(\neg\varphi)=n$).

A sentence $\varphi\in{\rm Sent}(\Sigma)$ is called {\em
$(P,\infty)$-{\rm (}co-{\rm )}rich} if ${\rm
RS}_P(\varphi)=\infty$ (respectively, ${\rm
RS}_P(\neg\varphi)=\infty$).

If $P=\mathcal{T}_\Sigma$ we write that $\varphi$ is
$(\alpha,n)$-{\rm (}co-{\rm )}rich instead of $(P,\alpha,n)$-{\rm
(}co-{\rm )}rich, and $\infty$-{\rm (}co-{\rm )}rich instead of
$(P,\infty)$-{\rm (}co-{\rm )}rich.

If for a property $P$ there is a $(P,\ast)$-{\rm (}co-{\rm )}rich
sentence $\varphi$, we say that $P$ has a {\em $(P,\ast)$-{\rm
(}co-{\rm )}rich sentence}, where $\ast=\alpha,n$ or
$\alpha=\infty$.

\medskip
By the definition if a sentence $\varphi$ is $(P,\alpha,n)$-rich
then ${\rm RS}(P_\varphi)=\alpha$, ${\rm ds}(P_\varphi)=n$.

\medskip
{\bf Theorem 3.4.} {\em $(1)$ If a property
$P\subseteq\mathcal{T}_\Sigma$ has a $(P,\alpha,m)$-rich sentence
$\varphi$ which is $(P,\beta,n)$-co-rich then ${\rm RS}(P)={\rm
max}\{\alpha,\beta\}$, ${\rm ds}(P)=m$ for $\alpha>\beta$, ${\rm
ds}(P)=n$ for $\alpha<\beta$, and ${\rm ds}(P)=m+n$ for
$\alpha=\beta$.

$(2)$ If for a property $P\subseteq\mathcal{T}_\Sigma$, ${\rm
RS}(P)=\alpha$ and ${\rm ds}(P)=n$, then for each sentence
$\varphi\in{\rm Sent}(\Sigma)$ the following assertions hold:

{\rm (i)} ${\rm RS}_P(\varphi)\leq\alpha$,

{\rm (ii)} if ${\rm RS}_P(\varphi)=\alpha$ then $\varphi$ is
$(P,\alpha,m)$-rich for some $m\leq n$, and for $m=n$ either
$\varphi\rhd_{\rm tt}P$ or $\varphi$ is $(P,\beta,k)$-co-rich for
some $\beta<\alpha$ and $k\in\omega$, and if $m<n$ then $\varphi$
is $(P,\alpha,n-m)$-co-rich.}

\medskip
Proof. (1) For the sentence $\varphi$ we have ${\rm
RS}_P(\varphi)=\alpha$, ${\rm ds}_P(\varphi)=m$, ${\rm
RS}_P(\neg\varphi)=\beta$, ${\rm ds}_P(\neg\varphi)=n$. It means
that $P$ is divided into two disjoint parts $P_\varphi$ and
$P_{\neg\varphi}$ with given characteristics ${\rm
RS}(P_\varphi)=\alpha$, ${\rm ds}(P_\varphi)=m$, ${\rm
RS}(P_{\neg\varphi})=\beta$, ${\rm ds}(P_{\neg\varphi})=n$.

If ${\rm RS}(P)=0$ then $|P|={\rm ds}(P)$, ${\rm
RS}_P(\varphi)={\rm RS}_P(\neg\varphi)=0$, ${\rm
ds}(P)=|P_\varphi|+|P_{\neg\varphi}|=m+n$.

If ${\rm RS}(P)>0$ then a tree witnessing the value ${\rm
RS}(P)=\gamma$ can be transformed step-by-step using theories
either in $P_\varphi$ or in $P_{\neg\varphi}$: in each step
witnessing ${\rm RS}(P)=\gamma$ there are infinitely many branches
of previous values related to $P_\varphi$ or to $P_{\neg\varphi}$.

In the first case, related to $P_\varphi$, we have $\gamma=\alpha$
and in the second case, related to $P_{\neg\varphi}$,
$\gamma=\beta$. If $\alpha>\beta$ a tree for $P_\varphi$ witnesses
$l=m$. If $\alpha<\beta$ a tree for $P_{\neg\varphi}$ witnesses
$l=n$. If $\alpha=\beta$ then both trees for $P_\varphi$ and for
$P_{\neg\varphi}$ witness $\gamma=\alpha$ and $l=m+n$, since there
are exactly $l+m$ $s$-definable subsets of $P$ having the rank
$\gamma$ and the degree $1$.

$(2)$ We have (i) by the monotony of rank (if $P_1\subseteq P_2$
then ${\rm RS}(P_1)\leq{\rm RS}(P_2)$) and the inclusion
$P_\varphi\subseteq P$. (ii) holds by the monotony of degree for a
fixed rank (if $P_1\subseteq P_2$ and ${\rm RS}(P_1)={\rm
RS}(P_2)\in{\rm Ord}$ then ${\rm ds}(P_1)\leq{\rm ds}(P_2)$) and
the inclusion $P_\varphi\subseteq P$. Besides if ${\rm
ds}(P_\varphi)=m=n={\rm ds}(P)$ then $P$ can not have a tree in
$P_{\neg\varphi}=P\setminus P_\varphi$ witnessing ${\rm
RS}(P_{\neg\varphi})=\alpha$ since otherwise ${\rm ds}(P)$ should
be more than $n$. Therefore either $\neg\varphi$ is
$P$-inconsistent, i.e., $\varphi\rhd_{\rm tt}P$, or $\varphi$ is
$(P,\beta,k)$-co-rich for some $\beta<\alpha$ and $k\in\omega$. If
$m<n$ then $\varphi$ is $(P,\alpha,n-m)$-co-rich in view of
(1).~$\Box$

\medskip
By Theorem 3.4 for any $e$-totally transcendental property $P$ and
any $\alpha\leq{\rm RS}(P)$ there are $s$-definable subfamilies
$P_\varphi$ with ${\rm RS}(P_\varphi)=\alpha$. Similarly all
values $m\leq{\rm ds}(P)$ are also realized by appropriate
$s$-definable subfamilies.

Thus the {\em spectrum} ${\rm Sp}_{\rm Rd}(P)$ for the pairs
$({\rm RS}_P(\varphi),{\rm ds}_P(\varphi))$ with nonempty
$P_\varphi$ forms the set
\begin{equation}\label{eqs1}\{({\rm RS}(P),m)\mid 1\leq m\leq{\rm
ds}(P)\}\cup\{(\alpha,m)\mid \alpha<{\rm
RS}(P),m\in\omega\setminus\{0\}\},\end{equation} which is an
initial segment $O[(\beta,n)]$ consisting of all pairs
$(\alpha,m)\in{\rm Ord}\times(\omega\setminus{0})$ with
$\alpha\leq\beta$ and $m\leq n$ for $\alpha=\beta$, ${\rm
RS}(P)=\beta$, ${\rm ds}(P)=n$.

\medskip
{\bf Remark 3.5.} If ${\rm RS}(P)=\infty$ then $s$-definable
subfamilies $P_\varphi$ can have only values ${\rm
RS}(P_\varphi)=\infty$ or both the value ${\rm RS}(P)=\infty$ and
pairs forming some initial segment $O[(\beta,n)]$.

Indeed, let $P$ be a family of theories in a language $\Sigma$ of
independent $0$-ary predicates $Q_l$, $l\in\lambda$,
$\lambda\geq\omega$, such that each sentence
$Q^{\delta_1}_{i_1}\wedge\ldots\wedge Q^{\delta_k}_{i_k}$,
$i_1<\ldots<i_k<\lambda$, $\delta_1,\ldots,\delta_k\in\{0,1\}$,
$k\in\omega$, is $P$-consistent. Each $P$-consistent sentence
$\varphi\in{\rm Sent}(\Sigma)$ is divided into $2$-tree witnessing
that ${\rm RS}_P(\varphi)=\infty$. In such a case we say that the
spectrum ${\rm Sp}_{\rm Rd}(P)$ equals $\{\infty\}$.

The family $P$ above can be extended by a family $P'$ with
dependent predicates $Q_l$ producing given ${\rm RS}$-rank and
${\rm ds}$-degree for a subfamily with, say, $Q_0\leftrightarrow
Q_1$. Therefore the arguments for Theorem 3.4 produce an initial
segment $O[(\beta,n)]$ for the spectrum ${\rm Sp}_{\rm Rd}(P')$ of
$s$-definable family $P'$. Thus, ${\rm Sp}_{\rm Rd}(P\cup
P')=O[(\beta,n)]\cup\{\infty\}$.

Since each nonempty $s$-definable subfamily has a spectrum of the
form $O[(\beta,n)]$, or $\{\infty\}$, or
$O[(\beta,n)]\cup\{\infty\}$, initial segments are well-ordered
with respect to the relation $\subseteq$, and the ordinal ${\rm
RS}(P)$-ranks are bounded by $|\mathcal{T}_\Sigma|\leq 2^{{\rm
max}\{|\Sigma|,\omega\}}$, all values ${\rm Sp}_{\rm Rd}(P)$, for
nonempty properties $P\subseteq\mathcal{T}_\Sigma$, are exhausted
by these three possibilities, and we obtain the following:

\medskip
{\bf Theorem 3.6.} {\em For any nonempty property
$P\subseteq\mathcal{T}_\Sigma$ one of the following possibilities
holds for some $\beta\in{\rm Ord}$ and $n\in\omega\setminus\{0\}$:

$(1)$ ${\rm Sp}_{\rm Rd}(P)=O[(\beta,n)]$,

$(2)$ ${\rm Sp}_{\rm Rd}(P)=\{\infty\}$,

$(3)$ ${\rm Sp}_{\rm Rd}(P)=O[(\beta,n)]\cup\{\infty\}$.

All possibilities above are realized by appropriate languages
$\Sigma$ and properties $P\subseteq\mathcal{T}_\Sigma$.}

\medskip
{\bf Remark 3.7.} Any value ${\rm Sp}_{\rm Rd}(P)$ can be
naturally extended till $\overline{\rm Sp}_{\rm Rd}(P)={\rm
Sp}_{\rm Rd}(P)\cup\{-1\}$ corresponding to the value ${\rm
RS}(\emptyset)=-1$ of the empty subfamily of $\mathcal{T}_\Sigma$.
It is also natural to put $\overline{\rm Sp}_{\rm Rd}(P)=\{-1\}$
for empty $P\subset\mathcal{T}_\Sigma$. In view of Theorem 3.6 we
have the following description of extended spectra $\overline{\rm
Sp}_{\rm Rd}(P)$:

i) $\overline{\rm Sp}_{\rm Rd}(P)=\{-1\}$,

ii) $\overline{\rm Sp}_{\rm Rd}(P)=O[(\beta,n)]\cup\{-1\}$,

iii) $\overline{\rm Sp}_{\rm Rd}(P)=\{-1,\infty\}$,

iv) $\overline{\rm Sp}_{\rm Rd}(P)=O[(\beta,n)]\cup\{-1,\infty\}$.

\medskip
{\bf Theorem 3.8} \cite{dst}. {\em Let $\mathcal{T}$ be a family
of a countable language $\Sigma$ and with ${\rm
RS}(\mathcal{T})=\infty$, $\alpha$ be a countable ordinal,
$n\in\omega\setminus\{0\}$. Then there is a $d_\infty$-definable
subfamily $\mathcal{T}^\ast\subset\mathcal{T}$ such that ${\rm
RS}(\mathcal{T}^\ast)=\alpha$ and ${\rm ds}(\mathcal{T}^\ast)=n$.}

Theorems 3.6 and 3.8 immediately imply:

\medskip
{\bf Corollary 3.9.} {\em Let $\mathcal{T}$ be a family of a
countable language $\Sigma$ and with ${\rm
RS}(\mathcal{T})=\infty$, $\alpha$ be a countable ordinal,
$n\in\omega\setminus\{0\}$. Then there is a $d_\infty$-definable
property $P\subset\mathcal{T}$ such that ${\rm Sp}_{\rm
Rd}(P)=O[(\alpha,n)]$.}

\section{Spectra for cardinalities of definable subproperties}

In this section we study some refinements of the relation
$\rhd_{\rm pt}$.

For a cardinality $\lambda\geq 1$, a sentence $\varphi\in{\rm
Sent}(\Sigma)$ and a property $P\subseteq\mathcal{T}_\Sigma$ we
write $\varphi\rhd^\lambda_{\rm pt} P$ if $\varphi$ satisfies
exactly $\lambda$ theories in $P$, i.e., $|P_\varphi|=\lambda$.

By the definition if $P\ne\emptyset$ and $\varphi\rhd_{\rm tt} P$
then $\varphi\rhd^{|P|}_{\rm pt} P$, and conversely
$\varphi\rhd^{|P|}_{\rm pt} P$ implies $\varphi\rhd_{\rm tt} P$
for finite $P$. For infinite $P$ the converse implication can
fail. Moreover, since infinite sets can be divided into two parts
of same cardinality, one can easily introduce an expansion $P'$ of
$P$ by a $0$-ary predicate $Q$ such that $Q\rhd^{|P'|}_{\rm pt}
P'$ and $\neg Q\rhd^{|P'|}_{\rm pt} P'$, implying that
$Q\not\rhd_{\rm tt} P'$.

For a property $P$ we denote by ${\rm Sp}_{\rm pt}(P)$ the set
$\{\lambda\mid\varphi\rhd^\lambda_{\rm pt} P$ for some sentence
$\varphi\}$. This set is called the {\em ${\rm pt}$-spectrum} of
$P$.

By the definition $|P|\in{\rm Sp}_{\rm pt}(P)$ for any nonempty
$P$, and $\lambda\leq |P|$ for any $\lambda\in{\rm Sp}_{\rm
pt}(P)$, i.e., ${\rm sup}{\rm Sp}_{\rm pt}(P)=||P|+1|$.

A natural question arises on a description of ${\rm pt}$-spectra.

This question is easily answered for finite $P$, since in such a
case all subsets of theories are separated as $s$-definable
singletons from their complements, and we obtain the following:

\medskip
{\bf Proposition 4.1.} {\em For any finite property
$P\subseteq\mathcal{T}_\Sigma$, ${\rm Sp}_{\rm
pt}(P)=(|P|+1)\setminus\{0\}$.}

\medskip
The following assertion generalizes Proposition 4.1 in terms of
isolated points due Theorem 1.12:

\medskip
{\bf Proposition 4.2.} {\em If $P$ has exactly
$n\in\omega\setminus\{0\}$ isolated points then ${\rm Sp}_{\rm
pt}(P)\cap\omega=(n+1)\setminus\{0\}$.}

\medskip
Proof. Let $T_1,\ldots,T_n$ are isolated points in $P$. If
$k\in{\rm Sp}_{\rm pt}(P)\cap\omega$, i.e., there is a sentence
$\varphi$ with $\varphi\rhd^k_{\rm pt} P$ then $P_\varphi$
consists of isolated points $T_{i_1},\ldots,T_{i_k}$ since
elements of the finite set $P_\varphi$ are separated as
$s$-definable singletons from their complements in $P_\varphi$.
Then $k\in (n+1)\setminus\{0\}$. Conversely, if
$k\in\omega\setminus{\rm Sp}_{\rm pt}(P)$ then $k>n$ since each
$k\leq n$ equals $|P_\varphi|$ for some sentence $\varphi$,
implying $k\in{\rm Sp}_{\rm pt}(P)$.~$\Box$

\medskip
{\bf Proposition 4.3.} {\em For any nonempty property
$P\subseteq\mathcal{T}_\Sigma$ either ${\rm Sp}_{\rm
pt}(P)\cap\omega$ equals an initial segment $(n+1)\setminus\{0\}$
for some $n\in\omega\setminus\{0\}$ or ${\rm Sp}_{\rm
pt}(P)\cap\omega=\omega\setminus\{0\}$.}

\medskip
Proof. If there are sentences $\varphi$ with finite $P_\varphi$
then either there is $P_\varphi$ with the greatest finite
cardinality, implying ${\rm Sp}_{\rm
pt}(P)\cap\omega=(n+1)\setminus\{0\}$ by the arguments of
Proposition 4.2, or the finite cardinalities $|P_\varphi|$ are
unbounded that means ${\rm Sp}_{\rm
pt}(P)\cap\omega=\omega\setminus\{0\}$.~$\Box$

\medskip
{\bf Proposition 4.4.} {\em For any infinite property
$P\subseteq\mathcal{T}_\Sigma$ there is a nonempty set $Y\subseteq
|P|$ of infinite cardinalities such that either there is
$n\in\omega\setminus\{0\}$ with ${\rm Sp}_{\rm
pt}(P)=Y\cup(n+1)\setminus\{0\}$, or ${\rm Sp}_{\rm
pt}(P)=Y\cup\omega\setminus\{0\}$, or ${\rm Sp}_{\rm pt}(P)=Y$.
All values $Y\cup(n+1)\setminus\{0\}$, $Y\cup\omega\setminus\{0\}$
and $Y$, for a nonempty set $Y$ of infinite cardinalities and
$n\in\omega\setminus\{0\}$, are realized as ${\rm Sp}_{\rm pt}(P)$
for appropriate properties $P$.}

\medskip
Proof. Since ${\rm sup}{\rm Sp}_{\rm pt}(P)\leq |P|$ and
Proposition 4.3 describes all possibilities for ${\rm Sp}_{\rm
pt}(P)\cap\omega$, it suffices, for a nonempty set $Y$ of infinite
cardinalities and $n\in\omega\setminus\{0\}$, to find a property
$P_1$ with ${\rm Sp}_{\rm pt}(P_1)=Y\cup(n+1)\setminus\{0\}$, a
property $P_2$ with ${\rm Sp}_{\rm
pt}(P_2)=Y\cup\omega\setminus\{0\}$, and a property $P_3$ with
${\rm Sp}_{\rm pt}(P_2)=Y$. For the property $P_1$ one can take a
finite $n$-element family $Y_n$, expanded by a $0$-ary predicate
$Q_0$ marking all theories in $Y_n$, and extend $Y_n$ by families
$Y_\lambda$, for each $\lambda\in Y$, of $\lambda$ theories of
$\lambda$ independent $0$-ary predicates $Q^\lambda_i$ expanded by
a $0$-ary predicate $Q_\lambda$ marking all theories in
$Y_\lambda$. Any $P_1$-consistent sentence $\varphi$ satisfies
either $k\leq n$ theories in $Y_n$ or $\lambda$ many theories in
$Z_\lambda$ and possibly $\mu$ many theories in $Y_\mu$ for
finitely many $\mu<\lambda$. It means that the cardinalities
$|(P_1)_\varphi|$ witness the equality ${\rm Sp}_{\rm
pt}(P_1)=Y\cup(n+1)\setminus\{0\}$.

For the property $P_2$ we repeat the process for $P_1$ replacing
the part $Y_n$ by an $e$-minimal family $Y^\lambda_0$ consisting
of some $\lambda\in Y$ theories all of which are marked by new
$0$-ary predicate $Q_0$. Realizing this process we obtain that
$s$-definable sets are finite, cofinite, or consists of $\lambda$
theories for $\lambda\in Y$. Thus, ${\rm Sp}_{\rm
pt}(P_1)=Y\cup\omega\setminus\{0\}$.

For the property $P_3$ we repeat the process for $P_1$ without the
part $Y_n$ obtaining ${\rm Sp}_{\rm pt}(P_1)=Y$.~$\Box$

\medskip
Summarizing Propositions 4.1--4.4 we obtain the following theorem
describing ${\rm pt}$-spectra of properties.

\medskip
{\bf Theorem 4.5.} {\em For any nonempty property
$P\subseteq\mathcal{T}_\Sigma$ one of the following conditions
holds:

$(1)$ ${\rm Sp}_{\rm pt}(P)=(n+1)\setminus\{0\}$ for some
$n\in\omega\setminus\{0\}$; it is satisfied iff $P$ is finite with
$|P|=n$;

$(2)$ ${\rm Sp}_{\rm pt}(P)=Y\cup(n+1)\setminus\{0\}$ for some
nonempty set $Y\subseteq |P|$ of infinite cardinalities and
$n\in\omega\setminus\{0\}$;

$(3)$ ${\rm Sp}_{\rm pt}(P)=Y\cup\omega\setminus\{0\}$ for some
nonempty set $Y\subseteq |P|$ of infinite cardinalities;

$(4)$ ${\rm Sp}_{\rm pt}(P)=Y$ for some nonempty set $Y\subseteq
|P|$ of infinite cardinalities.

All values $(n+1)\setminus\{0\}$, $Y\cup(n+1)\setminus\{0\}$,
$Y\cup\omega\setminus\{0\}$, and $Y$, for a nonempty set $Y$ of
infinite cardinalities and $n\in\omega\setminus\{0\}$, are
realized as ${\rm Sp}_{\rm pt}(P)$ for an appropriate property
$P$.}

\medskip
The following assertion shows that $Y$ in Theorem 4.5 is finite
for a property $P$ with ${\rm RS}(P)=1$.

\medskip
{\bf Proposition 4.6.} {\em If ${\rm RS}(P)=1$ then ${\rm Sp}_{\rm
pt}(P)=Y\cup\omega\setminus\{0\}$ for some finite nonempty set
$Y\subseteq |P|$ of infinite cardinalities with $|Y|\leq{\rm
ds}(P)$.}

\medskip
Proof. Since ${\rm RS}(P)=1$, $P$ is divided into $k={\rm ds}(P)$
infinite $s$-definable $e$-minimal parts $P_1,\ldots,P_k$. Each
part $P_i$ has only finite and cofinite $s$-definable subset
producing ${\rm Sp}_{\rm
pt}(P_i)=\{|P_i|\}\cup\omega\setminus\{0\}$. Since each
$s$-definable subset of $P$ is a Boolean combination of
$s$-definable subsets of $P_i$ and $|Z_i\cup Z_j|={\rm
max}\{|Z_i|,|Z_j|\}$ for infinite $s$-definable $Z_i\subseteq
P_i$, $Z_j\subseteq P_J$, $i,j\leq k$, we obtain ${\rm Sp}_{\rm
pt}(P)=\{|P_1|,\ldots,|P_k|\}\cup\omega\setminus\{0\}$.~$\Box$

\medskip
{\bf Remark 4.7.} Describing possibilities for ${\rm pt}$-spectra
${\rm Sp}_{\rm pt}(P)$ we admit that properties $P$ may be not
$E$-closed. If we assume that $P$ is infinite and $E$-closed then
we have two cases: either $P$ is $e$-totally transcendental with
the least generating set of a cardinality $\mu_1\leq{\rm
max}\{|\Sigma|,\omega\}$ and with $\mu_2\leq{\rm
max}\{|\Sigma|,\omega\}$ accumulation points, or ${\rm
RS}(P)=\infty$ with $|P|\geq 2^\omega$ by Theorem 1.2. In the
first case values for $\lambda\in{\rm Sp}_{\rm pt}(P)$ are
exhausted by all cardinalities in $\omega\setminus\{0\}$ and by
some infinite cardinalities $\leq{\rm max}\{|\Sigma|,\omega\}$. In
particular, for countable $\Sigma$, since ${\rm RS}(P)$ is a
countable ordinal, we have ${\rm Sp}_{\rm
pt}(P)=(\omega+1)\setminus\{0\}$. In the second case values for
$\lambda\in{\rm Sp}_{\rm pt}(P)$ are exhausted by cardinalities in
$\omega\setminus\{0\}$ or by cardinalities of some its initial
segment, depending on existence of the least generating set for
$P$, and by some infinite cardinalities $\leq{\rm
max}\{|\Sigma|,\omega\}$ and cardinalities $\geq 2^\omega$. In
particular, for countable $\Sigma$, ${\rm Sp}_{\rm pt}(P)$
includes $2^\omega$, and possibly $\omega$ depending on existence
of infinite totally transcendental definable part. Thus, in
Theorem 4.5, for $E$-closed $P$ some cases are not realized:
$\{\omega\}\cup(n+1)\setminus\{0\}$, $\{\omega\}$, $\{\omega_1\}$,
$\{\omega,\omega_1\}$, etc.

\medskip
By Remark 4.7 using Theorem 1.2 we have the following:

\medskip
{\bf Theorem 4.8.} {\em For any nonempty $E$-closed property
$P\subseteq\mathcal{T}_\Sigma$ with at most countable language
$\Sigma$ one of the following possibilities holds:

$(1)$ ${\rm Sp}_{\rm pt}(P)=(n+1)\setminus\{0\}$ for some
$n\in\omega\setminus\{0\}$, if $P$ is finite with $|P|=n$;

$(2)$ ${\rm Sp}_{\rm pt}(P)=\{2^\omega\}\cup(n+1)\setminus\{0\}$
for some $n\in\omega$, if $P$ is infinite and has $n$ isolated
points;

$(3)$ ${\rm Sp}_{\rm pt}(P)=(\omega+1)\setminus\{0\}$, if $P$ is
infinite and totally transcendental;

$(4)$ ${\rm Sp}_{\rm pt}(P)=\{\omega,
2^\omega\}\cup\omega\setminus\{0\}$, if $P$ has an infinite
totally transcendental definable subfamily but $P$ itself is not
totally transcendental;

$(5)$ ${\rm Sp}_{\rm pt}(P)=\{2^\omega\}\cup\omega\setminus\{0\}$,
if $P$ has infinitely many isolated points but does not have
infinite totally transcendental definable subfamilies.}

\medskip
{\bf Remark 4.9.} Possibilities in Theorem 4.7 give low bounds for
correspondent cases in uncountable languages.

\section{$P$-generic sentences and $P$-generic theories}

{\bf Definition.} (Cf. \cite{Poizat, Truss, TeZi}) For a property
$P\subseteq\mathcal{T}_\Sigma$ a sentence $\varphi\in{\rm
Sent}(\Sigma)$ is called {\em $P$-generic} if ${\rm
RS}_P(\varphi)={\rm RS}(P)$, and ${\rm ds}_P(\varphi)={\rm ds}(P)$
if ${\rm ds}(P)$ is defined.

If $P=\mathcal{T}_\Sigma$ then we omit $P$ and a $P$-generic
sentence is called {\em generic}.

\medskip
By the definition we have the following:

\medskip
{\bf Proposition 5.1.} {\em Any $P$-generic sentence $\varphi$ is
$(P,{\rm RS}(P),{\rm ds}(P))$-rich if ${\rm RS}(P)$ is an ordinal,
and $(P,\infty)$-rich if ${\rm RS}(P)=\infty$. And vice versa,
each $(P,{\rm RS}(P),{\rm ds}(P))$-rich sentence, for an ordinal
${\rm RS}(P)$, is $P$-generic, and each $(P,\infty)$-rich
sentence, for ${\rm RS}(P)=\infty$, is $P$-generic.}

\medskip
{\bf Proposition 5.2.} {\em If $\varphi\rhd_{\rm tt} P$ then
$\varphi$ is $P$-generic.}

\medskip
In view of Proposition 5.2 any property
$P\subseteq\mathcal{T}_\Sigma$ has $P$-generic sentences.

\medskip
{\bf Corollary 5.3.} {\em If a property
$P\subseteq\mathcal{T}_\Sigma$ is finite and $\varphi\in{\rm
Sent}(\Sigma)$ then $\varphi\rhd_{\rm tt} P$ iff $\varphi$ is
$P$-generic.}

\medskip
Proof. By Proposition 5.2 it suffices to show that if $\varphi$ is
$P$-generic then $\varphi\rhd_{\rm tt} P$. If $P=\emptyset$ then
both ${\rm RS}(P)={\rm RS}_P(\varphi)=-1$ and $\varphi\rhd_{\rm
tt} P$. If $P$ consists of $n\geq 1$ theories, then ${\rm
RS}(P)=0$, ${\rm ds}(P)=n$. Assuming that $\varphi$ is $P$-generic
we have ${\rm RS}(P_\varphi)=0$, ${\rm ds}(P_\varphi)=n$ implying
that $\varphi$ belongs to all $n$ theories in $P$ and
$\varphi\rhd_{\rm tt} P$.~$\Box$

\medskip
{\bf Remark 5.4.} In view of Corollary 5.3 the converse
implication for Proposition 5.2 holds iff $P$ is finite. Indeed,
if $P=\mathcal{T}_\Sigma$ for a countable language $\Sigma$ with
${\rm RS}(P)=\infty$, which is characterized in Theorem 1.2, then
we can construct a $2$-tree of sentences $\varphi_\delta$,
$\delta\in {\,}<\omega 2$ witnessing this value of ${\rm RS}(P)$.
It means that $P$ is divided into two disjoint definable parts
$P_{\varphi_0}$ and $P_{\varphi_1}$ with ${\rm
RS}(P_{\varphi_0})={\rm RS}(P_{\varphi_1})=\infty$. Thus,
$\varphi_0$ and $\varphi_1$ are generic whereas
$\varphi_0\not\!\!\rhd_{\rm tt} P$ and $\varphi_1\not\!\!\rhd_{\rm
tt} P$. Moreover, this effect works both for an arbitrary property
$P$ with ${\rm RS}(P)=\infty$ and for an arbitrary property $P$
with ${\rm RS}(P)\in{\rm Ord}$ and ${\rm RS}(P)\geq 1$. In the
latter case we can remove a branch in the tree witnessing the
values ${\rm RS}(P)$ and ${\rm ds}(P)$ just considering a sentence
$\varphi_1\wedge\neg\psi$, where $\varphi_1$ is a tautology and
$P_\psi$ is nonempty with ${\rm RS}(P_\psi)<{\rm RS}(P)$. In such
a case $\varphi_1\wedge\neg\psi$ is $P$-generic and
$(\varphi_1\wedge\neg\psi)\not\!\!\rhd_{\rm tt} P$.~$\Box$

\medskip
In view of Remark 5.4 we have the following:

\medskip
{\bf Proposition 5.5.} {\em For a property
$P\subseteq\mathcal{T}_\Sigma$ there is a $P$-generic sentence
$\varphi\in{\rm Sent}(\Sigma)$ with minimal/least $P_\varphi$ iff
$P$ is finite. If that $\varphi$ exists then $P_\varphi=P$.}

\medskip
By Proposition 5.5 for a property $P\subseteq\mathcal{T}_P$ with
${\rm RS}(P)\geq 1$, $P$-generic sentences $\varphi$ produce
infinite descending chains of $s$-definable subfamilies
$P_\varphi$.

\medskip
{\bf Proposition 5.6.} {\em If for a property
$P\subseteq\mathcal{T}_\Sigma$, ${\rm RS}(P)=\alpha\in{\rm Ord}$,
${\rm ds}(P)=1$, then for any sentence $\varphi\in{\rm
Sent}(\Sigma)$ either ${\rm RS}_P(\varphi)=\alpha$ and ${\rm
ds}_P(\varphi)=1$, or ${\rm RS}_P(\neg\varphi)=\alpha$ and ${\rm
ds}_P(\neg\varphi)=1$.}

\medskip
Proof. By the conjecture for $P$ and the monotony for pairs of
values of ${\rm RS}$ and ${\rm ds}$ we have ${\rm
RS}_P(\varphi)\leq\alpha$ and ${\rm RS}_P(\neg\varphi)\leq\alpha$
for any $\varphi\in{\rm Sent}(\Sigma)$, and if the ${\rm RS}$-rank
equals $\alpha$ then the ${\rm ds}$-degree equals $1$. We can not
have ${\rm RS}_P(\varphi)<\alpha$ and ${\rm
RS}_P(\neg\varphi)<\alpha$ for $\varphi\in{\rm Sent}(\Sigma)$ by
Theorem 3.4. Thus ${\rm RS}_P(\varphi)=\alpha$ and ${\rm
ds}_P(\varphi)=1$, or ${\rm RS}_P(\neg\varphi)=\alpha$ and ${\rm
ds}_P(\neg\varphi)=1$. The latter conditions can not be satisfied
simultaneously since otherwise ${\rm ds}(P)\geq 2$.~$\Box$

\medskip
Since any property $P\subseteq\mathcal{T}_\Sigma$ with ${\rm
RS}(P)=\alpha\in{\rm Ord}$ is represented as a disjoint union of
${\rm ds}(P)$ $s$-definable subfamilies $P_1,\ldots,P_{{\rm
ds}(P)}$ with ${\rm RS}=\alpha$ and ${\rm ds}=1$, Proposition 5.6
immediately implies:

\medskip
{\bf Corollary 5.7.} {\em For any property
$P\subseteq\mathcal{T}_\Sigma$ with ${\rm RS}(P)=\alpha\in{\rm
Ord}$ and any sentence $\varphi\in{\rm Sent}(\Sigma)$ either
$\varphi$ is $P$-generic or $\neg\varphi$ is $P$-generic, or, for
${\rm ds}(P)>1$ with non-$P$-generic $\varphi$ and $\neg\varphi$,
$\varphi$ is represented as a disjunction of $k$
$(P,\alpha,1)$-rich sentences and $\neg\varphi$ is represented as
a disjunction of $m$ $(P,\alpha,1)$-rich sentences such that
$k+m={\rm ds}(P)$, $k>0$, $m>0$.}

\medskip
{\bf Remark 5.8.} By Proposition 5.6 for any property
$P\subseteq\mathcal{T}_\Sigma$ with ${\rm RS}(P)=\alpha\in{\rm
Ord}$ and ${\rm ds}(P)=1$ there is unique ultrafilter $U_P$
consisting of $P$-generic sentences. By Proposition 5.5 this
ultrafilter is principal if and only if $P$ is finite, i.e., in
such a case it is a singleton. Anyway $U_P$ produces a theory
$T\in\mathcal{T}_\Sigma$ consisting of $P$-generic sentences only.
This theory $T$ is called {\em $P$-generic}.

If $P$ is infinite then $T$ belongs to the $E$-closure ${\rm
Cl}_E(P)$ \cite{RSrank, cl} of $P$ as unique element of
$\alpha$-th Cantor--Bendixson derivative of ${\rm Cl}_E(P)$, i.e.,
an element of ${\rm Cl}_E(P)$ having Cantor--Bendixson rank ${\rm
CB}({\rm Cl}_E(P))=\alpha$ \cite{RSrank}.

If ${\rm ds}(P)>1$ we can divide $P$ into ${\rm ds}(P)$
$s$-definable parts $P_i$ with ${\rm RS}(P_i)={\rm RS}(P)$ and
${\rm ds}(P_i)=1$, each of which has unique $P_i$-generic theory
$T_i$. The set $\{T_1,\ldots,T_{{\rm ds}(P)}\}$ is called the {\em
set of $P$-generic theories.}

Thus we have the following:

\medskip
{\bf Proposition 5.9.} {\em Each $e$-totally transcendental
property $P$ has finitely many, exactly ${\rm ds}(P)\geq 1$,
$P$-generic theories. These theories have Cantor--Bendixson rank
${\rm CB}({\rm Cl}_E(P))=\alpha={\rm RS}(P)$.}

\medskip
Now we extend the results above on generic sentences and theories
for properties $P$ with ${\rm RS}(P)=\infty$.

\medskip
{\bf Proposition 5.10.} {\em If for a property
$P\subseteq\mathcal{T}_\Sigma$, ${\rm RS}(P)=\infty$, then for any
sentence $\varphi\in{\rm Sent}(\Sigma)$ either ${\rm
RS}_P(\varphi)=\infty$ or ${\rm RS}_P(\neg\varphi)=\infty$.}

\medskip
Proof. Assume that for a sentence $\varphi\in{\rm Sent}(\Sigma)$
we have ${\rm RS}_P(\varphi)<\infty$ and ${\rm
RS}_P(\neg\varphi)<\infty$. We can suppose that $\varphi$ and
$\neg\varphi$ are both $P$-consistent. Then $\varphi$ is
$(P,\alpha,m)$-rich and $(P,\beta,n)$-co-rich for some
$\alpha,\beta\in{\rm Ord}$, $m,n\in\omega$. Applying Theorem 3.4
we obtain ${\rm RS}(P)\in{\rm Ord}$ contradicting ${\rm
RS}(P)=\infty$.~$\Box$

\medskip
Proposition 5.10 immediately implies:

\medskip
{\bf Corollary 5.11.} {\em For any property
$P\subseteq\mathcal{T}_\Sigma$ with ${\rm RS}(P)=\infty$ and any
sentence $\varphi\in{\rm Sent}(\Sigma)$ either $\varphi$ is
$P$-generic or $\neg\varphi$ is $P$-generic.}

\medskip
{\bf Remark 5.12.} By Proposition 5.10 for any property
$P\subseteq\mathcal{T}_\Sigma$ with ${\rm RS}(P)=\infty$ there is
ultrafilters $U_P$ consisting of $P$-generic sentences. Since the
condition ${\rm RS}(P)=\infty$ implies the existence of $2$-tree
of $\Sigma$-sentences there are at least continuum many these
ultrafilters, and by Theorem 1.2 there are exactly continuum many
ones for the language $\Sigma$ with $|\Sigma|\leq\omega$. Each
$U_P$ produces a theory $T\in\mathcal{T}_\Sigma$ consisting of
$P$-generic sentences only. This theory $T$ is called {\em
$P$-generic}.

The $P$-generic theories form the perfect kernel with respect to
Cantor--Bendixson derivatives of ${\rm Cl}_E(P)$, i.e., the set of
element of ${\rm Cl}_E(P)$ having Cantor--Bendixson rank ${\rm
CB}({\rm Cl}_E(P))=\infty$.

\medskip
Applying Theorem 3.6 and summarizing Remarks 5.8 and 5.12 we
obtain the following:

\medskip
{\bf Theorem 5.13.} {\em $(1)$ For any nonempty property
$P\subseteq\mathcal{T}_\Sigma$ there are ${\rm ds}(P)$ $P$-generic
theories if $P$ is totally transcendental, and at least continuum
many if $P$ is not totally transcendental. In the latter case
either all theories in $P$ are $P$-generic if ${\rm Sp}_{\rm
Rd}(P)=\{\infty\}$, or $P$ has at least $\beta\cdot\omega+n$
non-$P$-generic theories if ${\rm Sp}_{\rm
Rd}(P)=O[(\beta,n)]\cup\{\infty\}$.

$(2)$ The ${\rm CB}$-rank of each $P$-generic theory equals ${\rm
RS}(P)$.}

\medskip
{\bf Definition} \cite{at}. For a property
$P\subseteq\mathcal{T}_\Sigma$ a sentence $\varphi\in{\rm
Sent}(\Sigma)$ is called {\em $P$-complete} if $\varphi$ isolates
a unique theory $T$ in $P$, i.e., $P_\varphi$ is a singleton. In
such a case the theory $T\in P_\varphi$ is called {\em
$P$-finitely axiomatizable} (by the sentence $\varphi$).

\medskip
Since $P$-finitely axiomatizable theories are isolated points we
have the following:

\medskip
{\bf Proposition 5.14.} {\em For any nonempty property
$P\subseteq\mathcal{T}_\Sigma$ a $P$-finitely axiomatizable theory
$T$ is $P$-generic iff $P$ is finite.}

\section{Conclusion}

We studied links between formulas and properties, considered and
described characteristics for properties with respect to
satisfying formulas, their ranks and degrees. Possibilities of
spectra for ranks, degrees, and cardinalities of definable
properties are shown. Generic formulas and theories are introduced
and characterized.

Possibilities for ranks and degrees for formulas and theories with
respect to given properties are described. Their top values form
generic sentences and theories which are also described and
characterized.

There are many natural model-theoretic and other properties that
can be studied and described in this context. In this case, the
relations $\rhd_{\rm pt}$ and $\rhd_{\rm tt}$ are preserved under
natural closures. At the same time the operator of $E$-closure
does not preserve a series of natural model-theoretic properties.
For instance, there are families of strongly minimal whose
accumulation points have the strict order properties. It implies
that families of $\omega$-stabile, superstable and stable theories
can be not $E$-closed in these classes. Natural questions arise on
characteristics and characterizations of families in these classes
$\mathcal{T}$ whose $E$-closures are contained in $\mathcal{T}$.

\textbf{Sergey Vladimirovich Sudoplatov}, Doctor of Sciences
(Physics and Mathematics), docent; Leading researcher, Sobolev
Institute of Mathematics, Novosibirsk; Head of Chair, Novosibirsk
State Technical University; Professor, Novosibirsk State
University. \newline e-mail: sudoplat@math.nsc.ru


\begin{thebibliography}{99}
\bibitem{ErPa}{\em Ershov Yu.~L.} Mathematical
logic~/ Yu.~L.~Ershov, E.~A.~Palyutin.~--- Moscow~: FIZMATLIT,
2011.~--- 356~p. [in Russian]

\bibitem{RSrank} {\em Sudoplatov S.~V.} Ranks for families of theories
and their spectra~/ S.~V.~Sudoplatov // arXiv:1901.08464v1
[math.LO].~--- 2019.~--- 17~p.

\bibitem{MarkhSud} {\em Markhabatov N.~D.} Ranks for families of all theories
of given languages~/ N.~D.~Markhabatov, S.~V.~Sudoplatov //
Eurasian Mathematical Journal.~--- 2021 (to appear).
arXiv:1901.09903v1 [math.LO].~--- 2019.~--- 9~p.

\bibitem{dst} {\em Markhabatov N.~D.} Definable subfamilies of theories,
related calculi and ranks~/ N.~D.~Markhabatov, S.~V.~Sudoplatov //
Siberian Electronic Mathematical Reports.~--- 2020.~---
Vol.~17.~--- P.~700--714.

\bibitem{at} {\em Sudoplatov S.~V.} Approximations of
theories~/ S.~V.~Sudoplatov // Siberian Electronic Mathematical
Reports.~--- 2020.~--- Vol.~17.~--- P.~715--725.

\bibitem{cl} {\em Sudoplatov S.~V.} Closures and generating sets related
to combinations of structures~/ S.~V.~Sudoplatov // Bulletin of
Irkutsk State University. Series ``Mathematics''.~--- 2016.~---
Vol.~16.~--- P.~131--144.

\bibitem{hataf} {\em Sudoplatov S.~V.} Hierarchy of families of theories and their rank
characteristics~/ S.~V.~Sudoplatov // Bulletin of Irkutsk State
University. Series ``Mathematics''.~--- 2020.~--- Vol.~33.~---
P.~80--95.

\bibitem{Poizat} {\em Poizat B.} Groupes Stables~/ B.~Poizat.~--- Villeurbanne~: Nur Al-Mantiq Wal-Ma'rifah, 1987.

\bibitem{Truss} {\em Truss J.~K.} Generic Automorphisms of Homogeneous
Structures~/ J.~K.~Truss // Proceedings of the London Mathematical
Society.~--- 1992.~--- Vol.~65, No.~3.~--- P.~121--141.

\bibitem{TeZi} {\em Tent K.} A Course in Model Theory
(Lecture Notes in Logic. No. 40)~/ K.~Tent, M.~Ziegler.~---
Cambridge~: Cambridge University Press, 2012.~--- 258~p.

\end{thebibliography}
\end{document}